# Certain properties of contra-$T^*_{12}$-continuous functions

*Hadi J. Mustafa*

*Dept. of Math., F. of Math. & Comp. Sci.*

*University of Kufa*

*Najaf, Iraq*

drhadi.mustafa@gmail.com

*Layth M. Alabdulsada*

*Institute of Mathematics*

*University of Debrecen*

*H-4002 Debrecen, P.O. Box 400 Hungary*

layth.muhsin@science.unideb.hu

***Abstract**—The concept of contra function was introduced by Dontchev [2], in this work, we use the notion of $T^*_{12}$-open to study a new class of function called a contra-$T^*_{12}$-continuous function as generalization of contra- continuous.*

## I. INTRODUCTION

In 1996, Dontchev [2] introduced contra-continuous functions. In [10], the authors introduced the concept of almost contra-$T^*$-continuous function. In this paper, we introduce a new class of function called contra-$T^*_{12}$-continuous function where $T_1$, $T_2$ are operators associated with the topology τ on X. Throughout the paper, the space X and Y or (X, Y) and (Y, δ) stand for topological space, let A be a subset of X. the closure of A and the interior of A will be denoted by Cl(A) and Int(A), respectively.

## II. PRELIMINARIES

In this section, we recall the basic facts and definitions needed in this work.

**2.1 Definition:** A subset A of a space X is said to be:

**i)** Semi-open [6] if A ⊆ Cl (Int(A)),

**ii)** Pre-open [7] if A ⊆ Int (Cl (A)),

**iii)** b-open [1] if A ⊆ Cl (Int (A))∪ Int (Cl(A)).

The complement of semi-open (pre-open, b-open) is said to be semi-closed (pre-closed, b-closed). The family of all semi-open (pre-open, b-open, semi-closed, pre-closed, b-closed) subset of a space X is denoted by SO(X)(PO(X), BO(X), SC(X), PC(X), BC(X), respectively).

**2.2 Definition [4]:** A function f : X→Y is called semi-continuous ( pre-continuous, b-continuous) if for each x ∈ X and each open set V of Y containing f(x), there exists U ∈ SO(X) (U ∈ PO(X), U ∈ BO (X))  such that f(U) ⊆ V.

**2.3 Definition:** A function f : X → Y is called contra-continuous [2] (contra-semi continuous [4], contra-pre-continuous [3], contra-b-continuous [5]) if $f^{-1}(V)$ is closed (semi-closed, pre-closed, b-closed, resp.) in X for each open set V of Y.

### III. OPERATOR TOPOLOGICAL SPACES

**3.1 Definition [8]:** Let (X, τ) be a topological space and let T: p(X) → p(X)  be a function (where p(X) is the power set of X) we say that T is an operator associated with the topology τ on X if W ⊆ T(W)  (W ∈ τ) and the triple ( X, τ, T) is called an operator topological space.

**3.2 Definition [9]:** Let (X, τ, T) be an operator topological space, let A ⊆ X

**i)** A is called T-open if given x ∈ A, then there exists V ∈ τ there exists x ∈ V ⊆T (V) ⊆ A.

**ii)** A is called T*-open if A ⊆ T (A) (A is not necessarily open).

**3.3 Remarks:**

**i)**   Every T-open set is open.

**ii)** Every open set is T*-open, so we have the following implications:

T-open    →     open    →    T*-open

**iii)**   Let (X, τ) be a topological space define     T: p(X) → p(X) as follows: T (A) = Int Cl(A) then T is an operator associated with the topology τ on X and the triple (X, τ, T) is an operator topological space.

As an example, we can suppose X = R, τ =$t_u$= the usual topology on R, if

T (A) = Int Cl (A),

then the triple (R, $t_u$, T) is an operator topological space,

notice that Q ⊂ R satisfies Q ⊆ Int (Cl (Q)), so Q is a T*-open (pre-open) which is not open.

**3.3 Definition:** Let (X, τ) be a topological space and let $T_1$, $T_2$ be two operators associated with the topology τ on X then (X, τ, $T_1$, $T_2$) is called a bi operator topological space.

**3.4 Definition:** Let  (X, τ, $T_1$, $T_2$) be an  operator topological space and let  A ⊆ X, we say that A is a $T^*_{12}$-open if A ⊆ $T_1$(A) U $T_2$(A), the complement of $T^*_{12}$-open is called $T^*_{12}$-closed for example if:

$T_1(A)$ = Cl (Int(A)),

$T_2(A) = Int (Cl (A))$, Then

$A \subseteq Cl (Int (A)) \cup Int (A)$,

this is the definition of b-open set.

Notice that every $T^*_1$-open ($T^*_2$-open) is $T^*_{12}$-open because if A is $T^*_1$-open then $A \subseteq T_1(A) \subseteq T_1(A) \cup T_2(A)$, so A will be $T^*_{12}$-open.

## IV. CONTRA-$T^*_{12}$-CONTINUOUS FUNCTIONS

In this section, we obtain some properties of contra-$T^*_{12}$-continuous functions.

**4.1 Lemma [1]:** Let $(X, \tau)$ be a topological space then:

1) The intersection of an open set and a b-open set is a b-open set.

2) The union of any family of b-open sets is a b-open set.

Now, we generalize Lemma 4.1 as follows:

**4.2 Lemma:** Let $(X, \tau, T_1, T_2)$ be a bi operator topological space assume that

$T_1(W \cap B) = T_1(W) \cap T_1(B)$, $W \in \tau$, $B \subseteq X$,

$T_2(W \cap B) = T_2(W) \cap T_2(B)$, $W \in \tau$, $B \subseteq X$, therefore:

1) The intersection of an open set and a $T^*_{12}$-open set is $T^*_{12}$-open.
2) The union of any family $T^*_{12}$-open sets is a $T^*_{12}$-open set.

**Proof:**

1) Let $W \in X$ be an open set and let V be a $T^*_{12}$–open set we have to prove that $W \cap V$ is also a $T^*_{12}$–open set. Since W is open then:

$W \subseteq T_1(W)$ ... (1)

$W \subseteq T_2(W)$ ... (2)

Since V is a $T^*_{12}$–open then

$V \subseteq T_1(V) \cap T_2(V)$ ... (3)

$W \cap V \subseteq W \cap (T_1(V) \cap T_2(V))$

$= (W \cap T_1(V)) \cup (W \cap T_2(V))$

$\subseteq (T_1(W) \cap T_1(V)) \cup (T_2(W) \cap T_2(V))$

$= (T_1(W \cap V)) \cup (T_2(W \cap V))$

Then $W \cap V$ is $T^*_{12}$-open set.

2) Let $\mathcal{L} = \{w_\alpha \mid \alpha \in I\}$ be any family of $T^*_{12}$-open sets we must prove that $\bigcup_\alpha w\alpha$ is also a $T^*_{12}$-open

$w_α ⊆ T_1(w_α) ∪ T_2(w_α)$ for each $α ∈ I$

$\bigcup_α wα ⊆ \bigcup_α (T_1(w_α) ∪ T_2(w_α))$

$\quad\quad = \bigcup_α T_1(w_α) ∪ \bigcup_α T_2(w_α)$

Now $\bigcup_α T_1(w_α) = T_1(\bigcup_α w_α)$

Also $\bigcup_α T_2(w_α) = T_2(\bigcup_α w_α)$

Then $\bigcup_α w_α ⊆ T_1(\bigcup_α w_α) ∪ T_2(\bigcup_α w_α)$ and $\bigcup_α w_α$ is a $T^*_{12}$-open.

### 4.3 Remarks:

**i)** The intersection of two $T^*_{12}$-open is not necessarily $T^*_{12}$-open, so the collection of all $T^*_{12}$-open sets is not necessarily a topology on X.

Let $τ^*_{(12)}$ be the topology generated by the collection of all $T^*_{12}$-open sets.

**ii)** The intersection of any collection of $T^*_{12}$-closed sets is $T^*_{12}$-closed. Let $T^*_{12}$-Cl(B)-intersection of all $T^*_{12}$-closed sets containing B.

Recall that for a function f: $X → Y$, the subset $\{(x, f(x)) \mid x ∈ X\} ⊆ X × Y$ is called the graph of f and denoted by G (f).

**4.4 Definition:** Let $f:(X, τ, T_1, T_2) → (Y, δ)$ be a function the graph G(f) of f is said to be contra-$T^*_{12}$-closed graph if for each $(x, y) ∈ (X × Y)-G(f)$ there exists U which is $T^*_{12}$-open containing x and a closed set V of Y containing y such that $(U×V) ∩ G(f) = ∅$. The implies that $f(U) ∩ V=∅$.

**4.5 Definition**: A space X is said to be contra-compact if every closed cover of X has a finite sub cover.

**4.6 Theorem :** Let $(X, τ, T_1, T_2)$ be a bi operator topological space and suppose $f :(X, τ, T_1, T_2) → (Y, δ)$ has a contra-$T^*_{12}$-closed graph, then the inverse image of a contra–compact set A of Y is $T^*_{12}$-closed in X.

**Proof:** Assume that A is contra-compact set of A and $x ∉ f^{-1}(A)$ for each $a ∈ A$, $(x, a) ∉ G (f)$. Then there exists $U_a$ which is $T^*_{12}$-closed containing x and $V_a$ which is closed in Y containing a such that

$f (U_a) ∩ V_a = ∅$.

Consider $\mathcal{L} = \{A ∩ V_a \mid a ∈ A\}$ and $\mathcal{L}$ is a closed cover of the subspace A, but A is contra-compact then there exists $a_1, a_2, a_3...a_n$ such that

$A ⊆ \bigcup_{i=1}^{n} V_{ai}$.

Let $U = \bigcap_{i=1}^{n} U_{ai}$,

then U is $T^*_{12}$-closed containing x and $f (U) ∩ A = ∅$, therefore

$U ∩ f^{-1}(A) = ∅$ and hence $x ∉ T^*_{12}$-Cl $(f^{-1}(A))$, this show that $f^{-1}(A)$ is $T^*_{12}$-closed.

**4.7 Theorem :** Let Y be contra –compact space and let $( X, \tau^*_{(12)}, T_1, T_2)$ be operator topological space ,suppose $f : ( X, \tau^*_{(12)}, T_1, T_2) \to (Y, \delta)$ has a contra-$T^*_{12}$-closed graph then f is contra $T^*_{12}$-continuous.

**Proof:** First we show that an open set U of Y is contra –compact and let $\mathcal{L} = \{ V_\alpha \mid \alpha \in \Lambda \}$ be a cover of U by closed sets $V_\alpha$ of U for each $\alpha \in \Lambda$ , then there exists a closed set $K_\alpha$ of Y such that $V_\alpha = K_\alpha \cap U$ ,then the family $\{ K_\alpha \mid \alpha \in \Lambda \} \cup \{ U^c \}$ is closed cover of Y. But Y is contra-compact then there exists $\alpha_1, \alpha_2 \ldots \alpha_n$ such that

$Y = ( \bigcup_{i=1}^{n} K_{\alpha i}) \cup (U^c)$, hence

$U = \bigcup_{i=1}^{n} V_{\alpha i}$.

This show that U is contra-compact by (theorem 4.6) $f^{-1}(U)$ is a $T^*_{12}$-closed in X then for f is contra $T^*_{12}$-continuous.

**4.8 Theorem:** Let $f : ( X, \tau, T_1, T_2) \to (Y, \delta)$ be a function and $g : X \to X \times Y$ the graph function of f defined by $g(x) = (x, f(x))$ for every $x \in X$, if g is contra-$T^*_{12}$-continuous then f is contra-$T^*_{12}$-continuous.

**Proof:** Since g is contra-$T^*_{12}$-continuous then $f^{-1}(U) = g^{-1}(X \times U)$ is a $T^*_{12}$-closed in X. Then f is contra-$T^*_{12}$-continuous.

**4.9 Theorem :** If $f : ( X, \tau, T_1, T_2) \to (Y, \delta)$ is contra-$T^*_{12}$-continuous and $g : ( X, \tau, T_1, T_2) \to (Y, \delta)$ is contra-continuous and Y is Urysohn space then $E = \{ x \in X \mid f(x) = g(x) \}$ is $T^*_{12}$-closed in X.

**Proof**: Let $x \in E^c$, then $f(x) \neq g(x)$, since Y is a Urysohn then there exists open sets V and W such that $f(x) \in V$, $g(x) \in W$, and

$Cl(V) \cap Cl(W) = \emptyset$.

Since f is contra-$T^*_{12}$–continuous then $f^{-1}(Cl(V))$ is $T^*_{12}$-open in X and g is contra-continuous then $g^{-1}(Cl(W))$ is open in X, let $U = f^{-1}(Cl(V))$, $G = g^{-1}(Cl(W))$.

Then $x \in U \cap G = A$, where A is $T^*_{12}$-open in X and

$f(A) \cap g(A) \subseteq f(U) \cap g(G) \subseteq Cl(V) \cap Cl(W) = \emptyset$, hence

$f(A) \cap g(A) = \emptyset$ and $A \cap E = \emptyset$, $A \subseteq E^c$,

where A is $T^*_{12}$-open there for $x \notin T^*_{12}$-$Cl(E)$, then E is $T^*_{12}$-closed in X.

**4.10 Definition:** A subset A of operator topological space $(X, \tau, T_1, T_2)$ is said to be $T^*_{12}$-dense in X if $T^*_{12}$-$Cl (A) = X$.

**4.11 Remarks:** Let $(X, \tau)$ be a topological space define:

$T_1: p(X) \to p(X)$

$T_2: p(X) \to p(X)$ as follows

$T_1 (A) = Int (Cl (A))$

$T_2(A) = Cl(Int(A))$, then $T^*_{12}$-dense subset will be b-dense and $T^*_{12}$-Cl(A) will be b-Cl(A) so b-dense in X mean that b-Cl(A) = X.

**4.12 Corollary:** Let f: $(X, \tau, T_1, T_2) \to (Y, \delta)$ is contra-$T^*_{12}$-continuous and g: $(X, \tau, T_1, T_2) \to (Y, \delta)$ is contra continuous if Y is Urysohn and f = g on $T^*_{12}$-dense set A ⊆ X then f = g on X.

**Proof:** since f is contra -$T^*_{12}$-continuous and is contra continuous and Y is Urysohn by previous Theorem E = {x ∈ X: f(x) = g(x)} is a $T^*_{12}$-closed in X. We have f = g on $T^*_{12}$-dense set A ⊆ E, then X = $T^*_{12}$-Cl(A) ⊆ $T^*_{12}$-Cl(E) = E. Hence f = g on X.

**4.13 Definition:** A bi operator topological space $(X, \tau, T_1, T_2)$ is called $T^*_{12}$-connected if X is not the Union of two non-empty $T^*_{12}$-open sets.

**4.14 Theorem:** If f: $(X, \tau, T_1, T_2) \to (Y, \delta)$ is contra-$T^*_{12}$-continuous from a $T^*_{12}$-connected space onto Y, then Y is not a discrete space.

**Proof:** Suppose that *Y* is discrete. Let ∅ ≠ A ⊂ Y then A is proper nonempty open and closed subset of Y. Then $f^{-1}(A)$ is a proper nonempty $T^*_{12}$-clopen ($T^*_{12}$-open and $T^*_{12}$-closed) subset of X such that X = $f^{-1}(A) \cup (f^{-1}(A))^c$ which means that X is $T^*_{12}$-disconnected which is a contradiction. Hence Y is not discrete.

## REFERENCES

[1] D. Andrijevi´c, *"On b-open sets"*. Mat. Vesnik 48, 1996, 59-64.

[2] J. Dontchev, *"Contra-continuous functions and strongly S-closed spaces"*. Internat.J. Math. Math. Sci. 19, 1996, 303-310.

[3] S. Jafari and T. Noiri, *"On contra-precontinuous functions"*. Bull. Malaysian Math. Sc.Soc. 25 (2002), 115-128.

[4] E. Ekici and M. Caldas, "Slightly -continuous functions", Bol. Soc. Paran. Mat. 22(2) (2004), 63-74.

[5] A.A. Nasef, *"Some properties of contra-continuous functions"*. Chaos Solitons Fractals 24 (2005), 471-477.

[6] N. Levine, *"Semi-open sets and semi-continuity in topological spaces"*. Amer. Math. Monthly 70 (1963), 36-41.

[7] A. S. Mashhour, M. E. Abd El-Monsef and S. N. El-Deeb, *"On precontinuous and weak precontinuous functions"*. Proc. Math. Phys. Soc. Egypt 51 (1982), 47-53.

[8] Hadi J. Mustafa, A. Lafta. *"Operator topological space"*. Journal the college of Education, Al-Mustansiriya University. (2009)

[9] Hadi J. Mustafa, and A. Abdul Hassan, *T-open sets*. M.Sc thesis. Mu'ta University Jordan. (2004)

[10] Hadi J. Mustafa, and Layth M. Alabdulsada *"On Almost contra T*-continuous functions"*. Journal the college Mathematics and computer Sciences, Vol. 1, no. 6. (2012)